\documentclass[a4paper,DIV=12,10pt]{scrartcl}

\usepackage[utf8]{inputenc}
\usepackage{amsfonts,amssymb,amsmath}
\usepackage{mathtools}

\usepackage{epsfig}
\usepackage{MnSymbol}
\usepackage{booktabs}
\usepackage{cite}
\usepackage{float}
\usepackage{hyperref}
\usepackage{siunitx}

\usepackage{wrapfig}

\usepackage{graphics}
\usepackage{color}

\numberwithin{equation}{section}
\numberwithin{table}{section}
\numberwithin{figure}{section}

\renewcommand{\d}{\mathrm d}
\renewcommand{\vec}{\boldsymbol}

\newcommand{\R}{\mathbb{R}}

\usepackage{booktabs}

\usepackage{cleveref}

\usepackage{caption}
\usepackage{subcaption}

\begin{document}

\title{Efficiency of local Vanka smoother geometric multigrid preconditioning for space-time finite element methods to the Navier--Stokes equations}

\author{
	Mathias Anselmann$^\dag$, Markus Bause$^\ddag$\thanks{bause@hsu-hh.de (corresponding author)}\\
	{\small $^\dag,{}^\ddag$ Helmut Schmidt University, Faculty of
		Mechanical and Civil Engineering, Holstenhofweg 85,}\\ 
	{\small 22043 Hamburg, Germany}
}
\date{}

\maketitle

\begin{abstract}
Numerical simulation of incompressible viscous flow, in particular in three space dimensions, continues to remain a challenging task. Space-time finite element methods feature the natural construction of higher order discretization schemes. They offer the potential to achieve accurate results on computationally feasible grids. Linearizing the resulting algebraic problems by Newton's method yields linear systems with block matrices built of $(k+1)\times (k+1)$ saddle point systems, where $k$ denotes the polynomial order of the variational time discretization. We demonstrate numerically the efficiency of preconditioning GMRES iterations for solving these linear systems by a $V$-cycle geometric multigrid approach based on a local Vanka smoother. The studies are done for the two- and three-dimensional benchmark problem of flow around a cylinder. Here, the robustness of the solver with respect to the piecewise polynomial order $k$ in time is analyzed and proved numerically.   
\end{abstract}

\section{Introduction}

Let $\Omega \subset \R^d$ with $d=2,3$ denote a bounded Lipschitz domain and $I:=(0,T]$ for some final time $T>0$. We consider the numerical approximation by space-time finite element methods of solutions to the nonstationary Navier--Stokes system
\begin{equation}
	\label{Eq:NSE}
		\partial_t \vec{v} + (\vec{v} \cdot \nabla) \vec{v}  - \nu \Delta \vec{v} + \nabla p = \vec{f}\,, \quad 
		\nabla \cdot \vec{v} = 0
		\hspace*{3ex} \text{in } \Omega  \times T\,,
\end{equation}
equipped with the initial condition $\vec{v}(0) = \vec{v}_0$ in $\Omega$ and appropriate boundary conditions on the different parts of $\partial \Omega$, either the Dirichlet condition $\vec v = \vec g$ or the do-nothing condition $\vec T(\vec v,p)\cdot \vec n = \vec 0$ with the asymmetric stress tensor $\vec T(\vec v,p) = \nu \nabla \vec v - p \vec I$ and the outer unit normal vector $\vec n$. In Section~\ref{Sec:NumStudy}, the benchmark setting of flow around a cylinder (cf.~\cite{STDK96}) is chosen as $\Omega$. 

In the recent decade, the development and analysis of space-time finite element methods (STFEMs) have strongly attracted researchers' interest. For this we refer to, e.g., \cite{ABM17,AB20,AB21,AB22,HST14} and the references therein. STFEMs offer appreciable advantages like the natural construction of higher order schemes, the application of duality based concepts of a-posteriori error control and space-time mesh adaptivity \cite{BG10} and the natural discretization of coupled problems of multi-physics. Here, discontinuous Galerkin time discretizations and inf-sup stable finite element methods in space with discontinuous pressure approximation are employed. For the time variable, a Lagrangian basis with respect to the $(k+1)$ Gauss--Radau quadrature points of each subinterval is used. By the choice of a discontinuous temporal test basis, a time-marching scheme is obtained, with temporal degrees of freedom associated with the Gauss--Radau quadrature nodes of the subintervals. For details of the construction of the schemes we refer to \cite{AB22,AB21}. We linearize the resulting nonlinear algebraic problem by a damped version of Newton's method. This yields linear systems of equations with block matrices of $(k+1)\times (k+1)$ subsystems, with each of them having saddle point structure. Here, $k$ denotes the piecewise polynomial order of the time discretization; cf.~\cite{AB22,AB20,HST14}. The structure of each of the subsystems resembles the discretization of \eqref{Eq:NSE} by the implicit Euler method in time and inf-sup stable pairs of finite elements in space. We note that there exits a strong link between the implicit Euler and the lowest order discontinuous Gakerkin time discretization; cf.\ \cite{T96}. The solution of the problems demands for efficient and robust iteration schemes. 

Geometric multigrid (GMG) methods are known as the most efficient iterative techniques for the solution of large linear systems arising from the discretization of partial differential equations, including incompressible viscous flow. GMG methods are applied in many variants \cite{DJRWZ18,HKXZ16,JT00,T99}. They exploit different mesh levels of the underlying problem in order to reduce different frequencies of the error by employing a relatively cheap smoother on each grid level. Different iterative methods have been proposed in the literature as smoothing procedure. They range from low-cost methods like Richardson, Jacobi, and SOR applied to the normal equation of the linear system to  collective smoothers, that are based on the solution of small local problems. Here, we use a Vanka-type smoother \cite{V86,M06} of the family of collective methods. Nowadays, GMG methods are employed as preconditioner in Krylov subspace iterations, for instance GMRES iterations, to enhance the linear solver's robustness, which is also done here. Parallel implementations of GMG techniques on modern computer architectures show excellent scalability properties and their high efficiency has been recognized. Numerical evidence of these performance properties is presented in, e.g., \cite{AB22,GHJRW16,GRSWW15,J02,JT00}. Analyses of GMG methods (cf., e.g., \cite{DJRWZ18,HKXZ16,M06}) have been done for linear systems in saddle point form, with matrix $\vec{\mathcal A}= \begin{pmatrix} \vec A & \vec B^\top \\ \vec B & - \vec C \end{pmatrix}$ and symmetric and positive definite submatrices $\vec A$ and $\vec C$, arising for instance from mixed discretizations of the Stokes problem. If higher order discontinuous Galerkin time stepping is used, the resulting linear system matrix is a $(k+1)\times (k+1)$ block matrix with each block being of the form $\vec{\mathcal A}$. This imposes an additional facet of complexity on the geometric multigrid preconditioner. In this work, we study numerically the performance properties of GMRES iterations with GMG preconditioning for space-time finite element discretizations of \eqref{Eq:NSE} . In particular, higher order polynomial approximation of the temporal variable is investigated. The benchmark of flow around a cylinder \cite{STDK96} in two- and three space dimensions is chosen as a test problem. Beyond the expectable spatial grid independency of the number of iterations of the linear solver, that is already shown numerically in \cite{AB21,HST14}, the robustness of the solver with respect to the (piecewise) polynomial degree of the time discretization is analyzed and shown here.  

\section{Numerical scheme}
\label{Sec:Scheme}

For the time discretization, we decompose $I=(0,T]$ into $N$ subintervals $I_n=(t_{n-1},t_n]$, $n=1,\ldots,N$, where $0=t_0<t_1< \cdots < t_{N-1} < t_N = T$. For the space discretization, let $\{\mathcal{T}_l\}_{l=0}^{L}$ be the decomposition on every multigrid level of $\Omega$ into (open) quadrilaterals or hexahedrals, with $\mathcal T_l = \{K_i\mid i=1,\ldots , N^{\text{el}}_l\}$, for $l=0,\ldots,L$. The finest partition is $\mathcal T_h=\mathcal T_L$. We assume that all the partitions $\{\mathcal{T}_l\}_{l=0}^{L}$  are quasi-uniform with characteristic mesh size $h_l$ and $h_l=\gamma h_{l-1}$, $\gamma \in (0,1)$ and $h_0 = \mathcal O(1)$.  On the actual mesh level, we denote by $\vec V_h \times Q_h \subset \vec H^1(\Omega)\times L^2(\Omega)$ the finite element spaces that are built on the inf-sup stable pair of spaces $(\mathbb Q_r)^d \times \mathbb P_{r-1}^{\text{disc}}$ for some $r\geq 2$; cf.~\cite[p.~115]{J16}. By an abuse of notation, we skip the index $l$ of the mesh level when it is clear from the context. We let $R:=\operatorname{dim}\, \vec V_h$ and $S:=\operatorname{dim}\, Q_h$ as well as $\hat R = \operatorname{dim}((\mathbb Q_r)^d)$ and $\hat S = \operatorname{dim}(\mathbb P_{r-1}^{\text{disc}})$. For the application of Dirichlet boundary conditions we employ Nitsche's method. This is only done since our approach and its implementation are embedded in a more general framework for flow on evolving domains \cite{AB21}. We put
\begin{equation}
	\label{Def:B_GamD}
	\begin{aligned}
		B_{\Gamma_D}(\vec w,(\vec \psi_h,\xi_h)) : =   - \langle \vec w, \nu \nabla \vec \psi_h \cdot \vec n + \xi_h  \vec n  \rangle_{\Gamma_D}
		 + \gamma_1 \nu \langle h^{-1} \vec w , \vec 
		\psi_h  \rangle_{\Gamma_D}  + \gamma_2 \langle h^{-1} \vec w \cdot \vec n, \vec \psi_h \cdot \vec n 
		\rangle_{\Gamma_D} 
	\end{aligned}
\end{equation}
for $\vec w \in \vec H^{1/2}(\Gamma_D)$ and $(\vec \psi_h,\xi_h) \in \vec V_h \times Q_h$, where $\gamma_1>0$ and $\gamma_2> 0$ are numerical (tuning) parameters for the penalization, chosen as $\gamma_1 = \gamma_2 = 35$.  For $(\vec v_h,p_h)\in \vec V_h\times Q_h$ and $(\vec \psi_h,\xi_h) \in \vec V_h\times Q_h$, we define
\begin{align*}
	A_h((\vec v_h,p_h),(\vec \psi_h,\xi_h)) & \coloneqq 
	\begin{aligned}[t]
		&\langle (\vec v_h \cdot \vec \nabla) \vec v_h, \vec \psi_h \rangle
		+ \nu \langle \nabla \vec v_h , \nabla \vec \psi_h  \rangle
		-\langle p_h, \nabla \cdot \vec \psi_h \rangle
		+ \langle \vec \nabla \cdot \vec v_h, \xi_h \rangle \\
		& - \langle \nu \nabla \vec 
		v_h \cdot \vec n - p_h \vec n, \vec \psi_h\rangle_{\Gamma_D} + B_{\Gamma_D}(\vec v_h,(\vec \psi_h,\xi_h))\,,
		\end{aligned}\\
	L_h((\vec \psi_h,\xi_h);\vec f, \vec g)	& \coloneqq  \langle \vec f,  \vec \psi_h \rangle  + B_{\Gamma_D}(\vec g,(\vec \psi_h,\xi_h)) \,.
\end{align*}

Using discontinuous Galerkin time discretization \cite{AB21,HST14} with piecewise polynomials of order $k$ leads to finding, in each subinterval $I_n$, functions $(\vec v_{\tau,h},p_{\tau,h})\in \mathbb P_k(I_n;\vec V_{h}) \times \mathbb P_k(I_n;Q_h)$, such that for all $(\vec \psi_{\tau,h},\xi_{\tau,h}) \in \mathbb P_k(I;\vec V_h)\times \mathbb P_k(I_n;Q_h)$, 
	\begin{equation}
	\label{Eq_FullyDisc}
	\begin{aligned}
	\int_{t_{n-1}}^{t_{n}} \langle\partial_t \vec v_{\tau,h}, \vec \psi_{\tau,h} \rangle & +
	A_h((\vec v_{\tau,h},p_{\tau,h}),(\vec \psi_{\tau,h},\xi_{\tau,h})) \d t \\[1ex] & + \langle [\vec v_{\tau,h}]_{n-1}, \vec \psi _{\tau,h}(t_{n-1}^+) \rangle = \int_{t_{n-1}}^{t_{n}}
	L_h(\psi_{\tau,h}; \vec f, \vec g) \d t \,.
	\end{aligned}
\end{equation}
Here, $\mathbb P_k(I_n; B)$ denotes the space of all polynomials on $I_n$ of degree at most $k\in N_0$ with values in the Banach space $B$. The quantity $[\vec v_{\tau,h}]_{n-1}$ denotes, as usual, the jump of the discrete solution $\vec v_{\tau,h}$ at the left time point $t_{n-1}$ of $I_n$. 

\section{Newton-GMRES-GMG solver}
\label{Sec:Sol}

To derive the algebraic form of the variational equation \eqref{Eq_FullyDisc}, the discrete functions $\vec v_{\tau,h}$ and $p_{\tau,h}$ are firstly represented in a Lagrangian basis $\{\chi_{n,l}\}_{l=0}^k\subset \mathbb P_k(I_n;\R)$ with respect to the $(k+1)$ Gauss--Radau quadrature points of $I_n$ such that 
\begin{equation*}
\vec v_{\tau,h}{}_{|I_n}(\vec x,t) = \sum_{l=0}^k \vec v^{n,l}(\vec x) \chi_{n,l}(t) \quad \text{and}  \quad p_{\tau,h}{}_{|I_n}(\vec x,t) = \sum_{l=0}^k p^{n,l}(\vec x) \chi_{n,l}(t)\,.
\end{equation*}
The resulting coefficient functions $(\vec v^{n,l},p^{n,l})\in \vec V_h\times Q_h$, for $l=0,\ldots,k$, are then developed in terms of the finite element basis of $\vec V_h$ and $Q_h$, respectively. Letting $\vec V_h= \operatorname{span}\{\vec \psi_1,\ldots,\vec \psi_R \}$ and $Q_h=\operatorname{span}\{\xi_1,\ldots,\xi_S \}$, there holds that 
\begin{equation}
	\label{Eq:SpaceExp}
	\vec v^{n,l}(\vec x) = \sum_{r=1}^R v_{r}^{n,l} \vec \psi_r(\vec x) \quad \text{and} \quad p^{n,l}(\vec x) = \sum_{s=1}^S p_{s}^{n,l} \xi_s(\vec x)\,.
\end{equation}
With a suitably defined function $\vec F_n: \R^{(k+1)\cdot (R+S)} \rightarrow \R^{(k+1)\cdot (R+S)}$, the variational problem \eqref{Eq_FullyDisc} is recovered as 
\begin{equation}
\label{Eq:NLP}
\vec F_n (\vec X_n) = \vec 0 \,, 	\qquad \text{with solution vector} \quad \vec X_n^\top = (\vec v^{n,0}, \vec p^{n,0},\ldots,\vec v^{n,k}, \vec p^{n,k})
\end{equation}
and row subvectors $\vec v^{n,l}= (v_{1}^{n,l},\ldots,v_{R}^{n,l})$ and $\vec p^{n,l}= (p_{1}^{n,l},\ldots,p_{S}^{n,l})$, for $l=0,\ldots,k$, built from the coefficients of the expansions in \eqref{Eq:SpaceExp}. For the explicit representation of $\vec F_n$ in the case of a similar time discretization scheme, we refer to \cite{AB20}. We note that $\vec X^n$ comprises the (spatial) degrees of freedom for all $(k+1)$ Gauss--Radau nodes, representing the degrees of freedom in time, of the subinterval $I_n$. Thus, the approximations at these time points are computed simultaneously in \eqref{Eq:NLP}.

To solve the nonlinear problem \eqref{Eq:NLP}, a damped version of Newton's method is applied. The Jacobian matrix $\vec J_n(\vec X^m_n)$ of the resulting linear system $\vec J_n(\vec X^m_n) \vec C^m_n = - \vec F(\vec X^m_n)$, with iteration index $m$ of Newton's method, for the Newton correction $\vec C^m_n$ to the vector $\vec X^m_n$ of degrees of freedom in $I_n$, such that $\vec X_n^{m+1}= \vec X_n^m + \vec C_n^m$, exhibits a complex block structure 
\begin{equation}
\label{Eq:Jac}
\vec J_n(\vec X^m_n) = (\vec J^{n,m}_{a,b})_{a,b=0}^k \;\text{with submatrices $\vec J^{n,m}_{a,b}$ of saddle point form} \; \vec J^{n,m}_{a,b} = \begin{pmatrix} \vec D^{n,m}_{a,b}& \vec B^{n,m\,\top}_{a,b}\\[1ex]  \vec B^{n,m}_{a,b}& \vec 0  \end{pmatrix}\,.
\end{equation}
The explicit representation of the submatrices $\vec D_{a,b}^{n,m}$, for $a,b\in \{0,\ldots,k\}$, reads as   
\begin{align*}
	\large (\vec{D}^{n,m}_{a,b}\large)_{i,j} & :=
	\int_{t_{n-1}}^{t_n}
	\langle \partial_t \chi_{n,a} \vec \psi_{i}, \chi_{n,b} \vec \psi_{j}  \rangle
	+
	\langle (\vec v_{\tau,h}^m \cdot \vec \nabla) \chi_{n,a} \vec \psi_{i}, \chi_{n,b} \vec \psi_{j} \rangle
	+
	\langle (\chi_{n,a} \vec \psi_{i} \cdot \vec \nabla) \vec v_{\tau,h}^{m}, \chi_{n,b} \vec \psi_{j} \rangle\\[1ex]
	& \quad +
	\nu \langle \nabla \chi_{n,a} \vec \psi_{i} , \nabla \chi_{n,b} \vec \psi_{j}  \rangle
     + \langle \nu \chi_{n,a} \nabla \vec \psi_{i} \cdot \vec n, \chi_{n,b} \vec \psi_{j} \rangle_{\Gamma_D}
	-
	\langle \chi_{n,a} \vec \psi_{i}, \nu \chi_{n,b} \nabla \vec \psi_{j} \cdot \vec n \rangle_{\Gamma_D}
	\\[1ex]
	& \quad +
	\gamma_1 h^{-1} \nu \langle \chi_{n,a} \vec \psi_{i} , \chi_{n,b} \vec \psi_{j}  \rangle_{\Gamma_D} + 
	\gamma_2 h^{-1} \langle \chi_{n,a} \vec \psi_{i} \cdot \vec n, \vec \chi_{n,b} \vec \psi_{j} \cdot \vec n 
	\rangle_{\Gamma_D}
	\\[1ex]
	& \quad +
	\bigl< \chi_{n,a}(t_{n-1}^+) \vec \psi_{i},  \chi_{n,b}(t_{n-1}^+) \vec \psi_{j} \bigr>\,, \quad \text{for} \;\; i,j \in \{1,\ldots,R\} \,.
\end{align*}
The matrices $\vec{B}_{a,b}^{n,m\,\top}$, for $a,b\in \{0,\ldots,k\}$, represent the pressure-velocity coupling and are defined by 
\begin{equation*}
	\large( \vec{B}_{a,b}^{n,m\, \top}\large)_{i,j} : =
	\int_{t_{n-1}}^{t_n}
	\langle \chi_{n,a} \nabla \cdot \vec \psi_{i}, \chi_{n,b} \xi_j \rangle
	-
	\langle \chi_{n,a} \nabla \vec \psi_{i}, \chi_{n,b} \xi_j  \vec n  \rangle_{\Gamma_D}
	\,\d t
	\,, \; \text{for} \;\; i \in [1,R]\,, \; j \in [1,S] \,.
\end{equation*}

The linear system $\vec J_n(\vec X^n_m) \vec C^n_m= - \vec F(\vec X^n_m)$ is solved by GMRES iterations \cite{S03} that are preconditioned by a $V$-cycle geometric multigrid algorithm. For details of its parallel implementation we refer to \cite{AB22}. On the coarsest spatial mesh partition $\mathcal T_0$ of the multigrid hierarchy $\{\mathcal{T}_l\}_{l=0}^{L}$, a parallel direct solver is applied to the linear system $\vec J_0 \vec d_0 = \vec b_0$ for the defect vector $\vec d_0$. On the mesh partition $\mathcal T_l$,  for $l=1,\ldots, L$, with linear system $\vec J_l \vec d_l = \vec b_l$ for the defect vector $\vec d_l$, we define an elementwise Vanka operator  $S_K:\R^{(k+1)\cdot (R+S)} \rightarrow  \R^{(k+1)\cdot (\hat R + \hat S)}$ for $K\in \mathcal T_l$ by (cf., e.g., \cite{HST14,J02,V86})
\begin{equation}
\label{Eq:VS}
\vec S_K(\vec{d}) \coloneqq \vec R_K \vec d + \omega \vec J_K^{-1} \vec R_K  (\vec b_l-\vec J_l \vec d)\,,
\end{equation}
 with some underrelaxation factor $\omega$ chosen here as $\omega = 0.7$.  In \eqref{Eq:VS}, the $K$-local restriction operator $\vec R_K:\R^{(k+1)\cdot (R+S)} \rightarrow  \R^{(k+1)\cdot (\hat R + \hat S)}$ is defined by $(\vec R_K \vec d)_{\hat \mu} := \vec d_{\mathrm{dof}(K,\hat \mu)}$ for all $\hat \mu \in \hat Z :=\{1,\ldots,(k+1)\cdot(\hat R + \hat S)\}$ and assigns to a global defect vector $\vec d\in \R^{(k+1)\cdot (R+S)} $ the local block vector $\vec R_K \vec d\in \R^{(k+1)\cdot (\hat R + \hat S)}$ that contains all components of $\vec d$ that are associated with all degrees of freedom (for all $k+1$ Gauss--Radau points of $I_n$) belonging to the element $K$; cf.\ \eqref{Eq:NLP}. Here, the mapping $\mathrm{dof}:\mathcal T_l\times \hat Z \rightarrow Z_l$, $\mu = \mathrm{dof}(K,\hat \mu)\in Z_l(K)$, yields for a given element $K\in \mathcal T_l$ and a local degree of freedom with number $\hat \mu \in \hat Z$ the uniquely defined global number $\mu \in Z_l$, where $Z_l$ denotes the index set of all global degrees of freedom of $\mathcal T_l$ and $Z_h(K)$ is its subset of all global degrees of freedom belonging to the element $K\in \mathcal T_l$ (all velocity and pressure values for the $(k+1)$ Gauss--Radau time points). Further, we put $(\vec J_K)_{{\nu},{\mu}}:= (\vec J)_{\mathrm{dof}(K,\hat \nu,k),\mathrm{dof}(K,\hat \mu)}$ for all $\hat \nu,\hat \mu \in \hat Z$. For the computation of the inverse of the local Jacobian $\vec J_K^{-1}$ we use LAPACK routines. The $J_{\max}$ pre- and post-smoothing iterations by \eqref{Eq:VS}, where $J_{\max}=4$ is chosen here, are done on each multigrid level $\mathcal T_l$, for $l=1,\ldots, L$, with $\vec J \vec d = \vec b$ and read as ($\vec d^0= \vec 0$)
\begin{subequations}
\begin{alignat}{2}
\label{Eq:S0}
&\textbf{for}\;\; j=1,\ldots, J_{\max}\;\; {\color{blue} \{}\textbf{for}\;\; \mu = 1,\ldots, N^{\text{dofs}}\;\; \{\vec p_{\mu} = 0;\; \vec z_\mu=0;\}\\
\label{Eq:S1}
&\quad\textbf{for}\;\; i=1,\ldots, N^{\text{el}}\;\; {\color{red}\{}\vec y_{K_i} = \vec S_{K_i}(\vec{d}^{j-1});\\
\label{Eq:S2}
&\quad \quad \textbf{for}\;\; \hat \mu = 1,\ldots, (k+1)\cdot (\hat R +\hat S)\;\; \{\vec z_{\mathrm{dof}(K_i,\hat \mu)} \texttt{+=} \vec (\vec E_{K_i} \vec y_{K_i})_{\mathrm{dof}(K_i,\hat \mu)}; \; \vec p_{\mathrm{dof}(K_i,\hat \mu)}\texttt{++};\}{\color{red}\}}\\
\label{Eq:S3}
& \quad \textbf{for}\;\; \mu = 1,\ldots, N^{\text{dofs}}\;\;  \{\vec d_\mu^j = \vec z_\mu/\vec p_\mu;\}{\color{blue}\}}
\end{alignat}
\end{subequations}
In \eqref{Eq:S1}, the local Vanka smoother \eqref{Eq:VS} is applied to the former defect $\vec d^{j-1}$. In  \eqref{Eq:S2} the local update $\vec y_{K_i}$ is assigned to a global auxiliary vector $\vec z$ and the vector $\vec p$ counting the updates of the degrees of freedom of $K_i$ is incremented. Finally, in \eqref{Eq:S3} the arithmetic mean of the local updates of the degrees of freedom is assigned to the update of the global defect vector $\vec d^j$. In \eqref{Eq:S2}, the $K$-dependent extension operator $\vec E_K:\R^{(k+1)\cdot (\hat R + \hat S)} \rightarrow \R^{(k+1)\cdot (R+S)}$ is defined by $$(\vec E_K \vec y_K)_\mu = \left\{\begin{array}{@{}ll} (\vec y_K)_{\hat \mu}\,, & \text{if}\; \exists \hat \mu \in \hat Z:\; \mu = \operatorname{dof}(K,\hat \mu)\,, \\ 0\,, & \text{if}\; \mu \not\in Z_l(K) \,.\end{array} \right.$$  
We note that the averaging of \eqref{Eq:S2} and \eqref{Eq:S3}, used instead of overwriting successively the degrees of freedom on the element faces within the element loop of \eqref{Eq:S1} and \eqref{Eq:S2}, improves the performance properties of the GMRES--GMG linear solver.

\section{Numerical study}
\label{Sec:NumStudy}

Here, we analyse the performance properties of the proposed algorithm. This is done for the benchmark problem \cite{STDK96} of flow around a cylinder in two and three space dimensions and an increasing polynomial order $k$ in time, impacting strongly the complexity of the Jacobian matrix $\vec J$ in \eqref{Eq:Jac}. In \cite{AB22}, it is shown numerically for the dG($1$) time discretization, corresponding to $k=1$ in \eqref{Eq_FullyDisc}, that the performance properties of the GMRES--GMG solver are independent from the resolution and refinement of the spatial grid. Here, we demonstrate that even a robustness of the solver with respect to the temporal approximation order is expectable. In all the numerical experiments, the stopping criterion of the Newton iteration is given by the condition of an absolute residual smaller than \num{1e-8}, for the GMRES iterations an absolute residual smaller than \num{1e-9} is prescribed.

\subsection{2D flow benchmark}
\label{sec:NumStudy_2D}

Firstly, we consider the two-dimensional flow benchmark of \cite{STDK96} in the \emph{2D-2} unsteady setting, that is characterized by the Reynolds number $Re=100$. Here, we focus on evaluating the performance properties of the algebraic solver of Section~\ref{Sec:Sol}. The polynomial degree $k$ of the approximation space in \eqref{Eq_FullyDisc} is increased from $k=0$ to $k=4$. The final simulation time is chosen to $T=7$ and the time step size is $\tau_n = 0.1$. For the approximation in space, we put $r=4$, such that the element pair $(\mathbb Q_4)^2\times \mathbb P_3^{\text{disc}}$ is used. For $T=7$ the flow profile is not fully developed yet. For this reason we do not report any values for the drag and lift coefficient as the goal quantities. For values of the drag and lift coefficient, computed by our approach in the case $k=1$, we refer to \cite{AB22}. \Cref{tab:2d_DFG_results} summarizes the convergence statistics of the Newton-GMRES-GMG solver. We nicely observe the solver's robustness with respect to the polynomial degree of the time discretization and, thereby, with respect to the $(k+1)\times (k+1)$ block structure of the Jacobian matrix in \eqref{Eq:Jac}. We note that the GMG preconditioned GMRES solver only needs, in the mean, 2 or 3 iterations for convergence and, thereby, is highly efficient.   

{
	\sisetup{scientific-notation = false,
		round-mode=places,
		round-precision=4,
		output-exponent-marker=\ensuremath{\mathrm{e}},
		table-figures-integer=1, 
		table-figures-decimal=3, 
		table-figures-exponent=1, 
		table-sign-mantissa = false, 
		table-sign-exponent = true, 
		table-number-alignment=center} 
	
	\begin{table}[h!t]
		\caption{Piecewise polynomial order $k$ in time, space-time degrees of freedom on each time interval $I_n$, average number $\bar{n}_{\text{Newton}}$ of Newton steps per time step and average number $\bar{n}_{\text{GMRES}}$ of GMRES iterations per Newton step for the entire time interval $t\in (0,T]$.}
		\centering
		\small
		\begin{tabular}{c@{\hskip 4ex} c@{\hskip 4ex} c@{\hskip 4ex} c@{\hskip 4ex} c}
	\toprule
	{$k$} & {DoFs$_{I_n}$} & {$\bar{n}_{\text{Newton}}$} & {$\bar{n}_{\text{GMRES}}$}  \\
	%
	\midrule
	{0} & 993024  & 1.28 & 2 \\
	{1} & 1986048 & 3.39 & 3 \\
	{2} & 2979072 & 3.68 & 3 \\
	{3} & 3972096 & 3.68 & 3 \\
	{4} & 4965120 & 3.62 & 3 \\
	\bottomrule
\end{tabular}
		\label{tab:2d_DFG_results}
	\end{table}
}

\subsection{3D flow benchmark}
\label{sec:NumStudy_3D}

Next, we consider the three-dimensional flow benchmark of \cite{STDK96} in the \emph{3D-2Z} setting, characterized again by the Reynolds number $Re = 100$. For the spatial discretization we put $r=2$, corresponding to the finite element pair $(\mathbb Q_2)^3\times \mathbb P_1^{\text{disc}}$. We let $T = 8.5$ and adapt the time step size $\tau_n$ on $I=(0,T]$. On the subinterval $I_i=(0,7]$ we put $\tau_n = 0.1$, on $I_m=(7,8]$ we set $\tau_n = 0.05$ and, finally, on $I_e=(8,8.5]$ we set $\tau_n = 0.01$. This adaptation is motivated by the idea, that a coarser time mesh is used during the initial period when the flow develops and a very fine time mesh is chosen in the final part of the simulation, when the goal quantities of interest are computed and monitored. The initial velocity is set to $\vec{v}(0) = \vec{0}$ and the inflow boundary condition is prescribed for $t >0 $, without any (smooth) startup  process for the application of the boundary condition. For the computation of the drag and lift coefficient, the surface integrals defining the drag and lift forces are approximated. Further, we recall that Nitsche's method, based on adding the terms defined in \eqref{Def:B_GamD} to the variational problem, is employed for the implementation of Dirichlet boundary conditions. 

\begin{figure}[!htb]
	\centering
	\subcaptionbox{\centering Drag coefficient in the initial phase. \label{fig:drag_start}}
	[0.47\columnwidth]
	{\includegraphics[width=0.47\columnwidth,keepaspectratio]{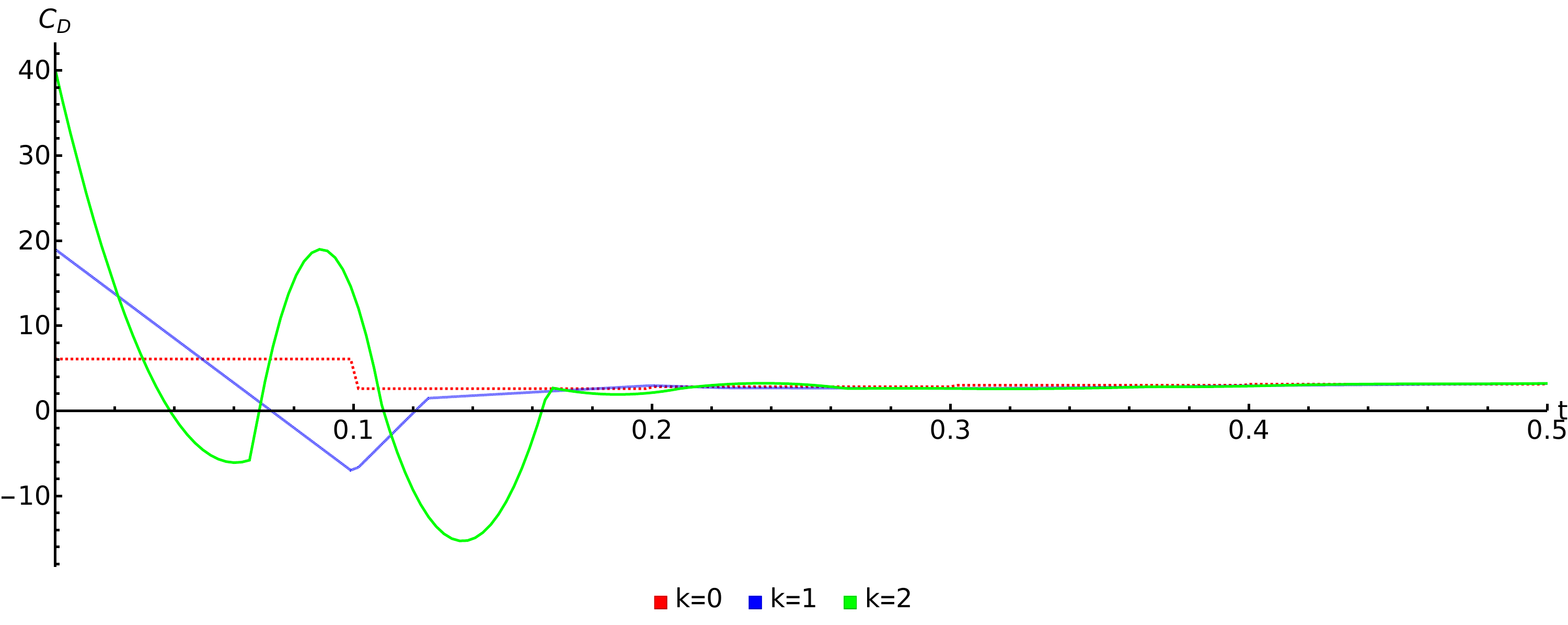}
	}
	\quad
	\subcaptionbox{\centering Lift coefficient in the initial phase. \label{fig:lift_start}}
	[0.47\columnwidth]
	{\includegraphics[width=0.47\columnwidth,keepaspectratio]{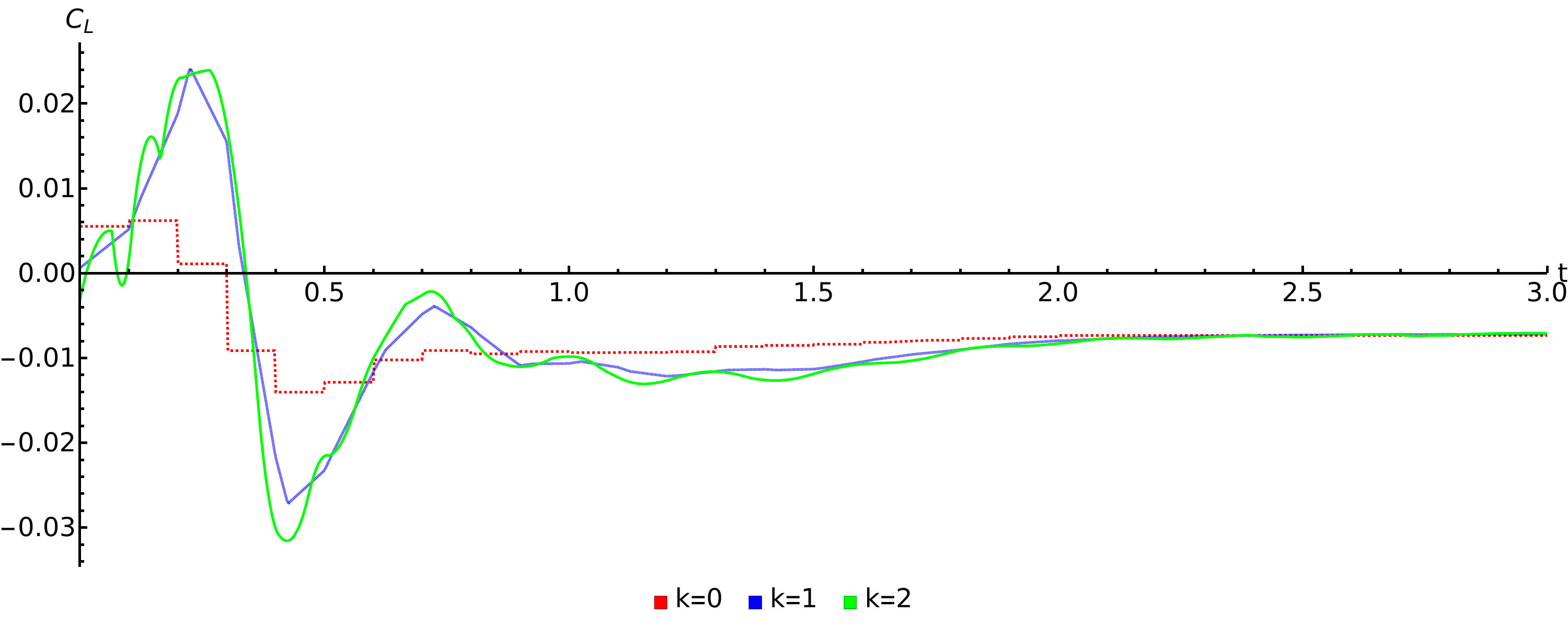}
	}
	\\[2ex]
	\subcaptionbox{\centering Drag coefficient for fully developed flow. \label{fig:drag_end}}
	[0.47\columnwidth]
	{\includegraphics[width=0.47\columnwidth,keepaspectratio]{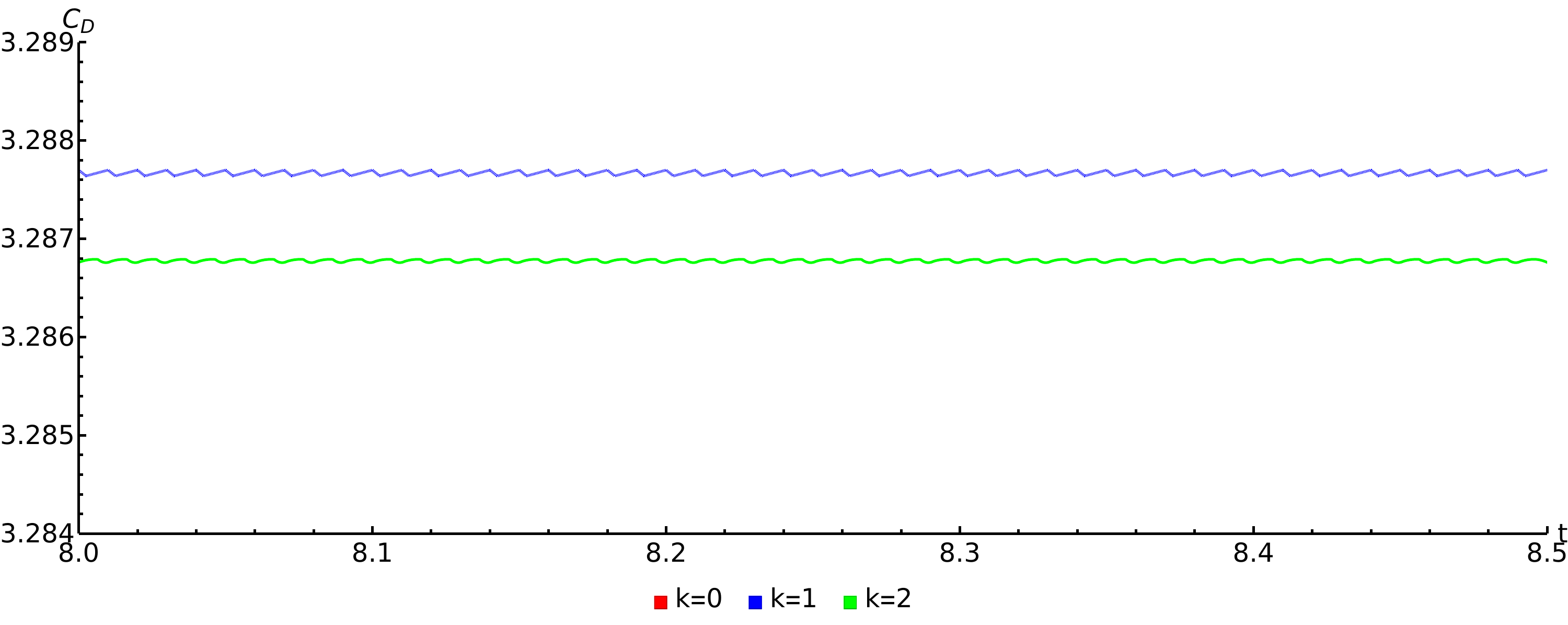}
	}
	\quad
	\subcaptionbox{\centering Lift coefficient for fully developed flow. \label{fig:lift_end}}
	[0.47\columnwidth]
	{\includegraphics[width=0.47\columnwidth,keepaspectratio]{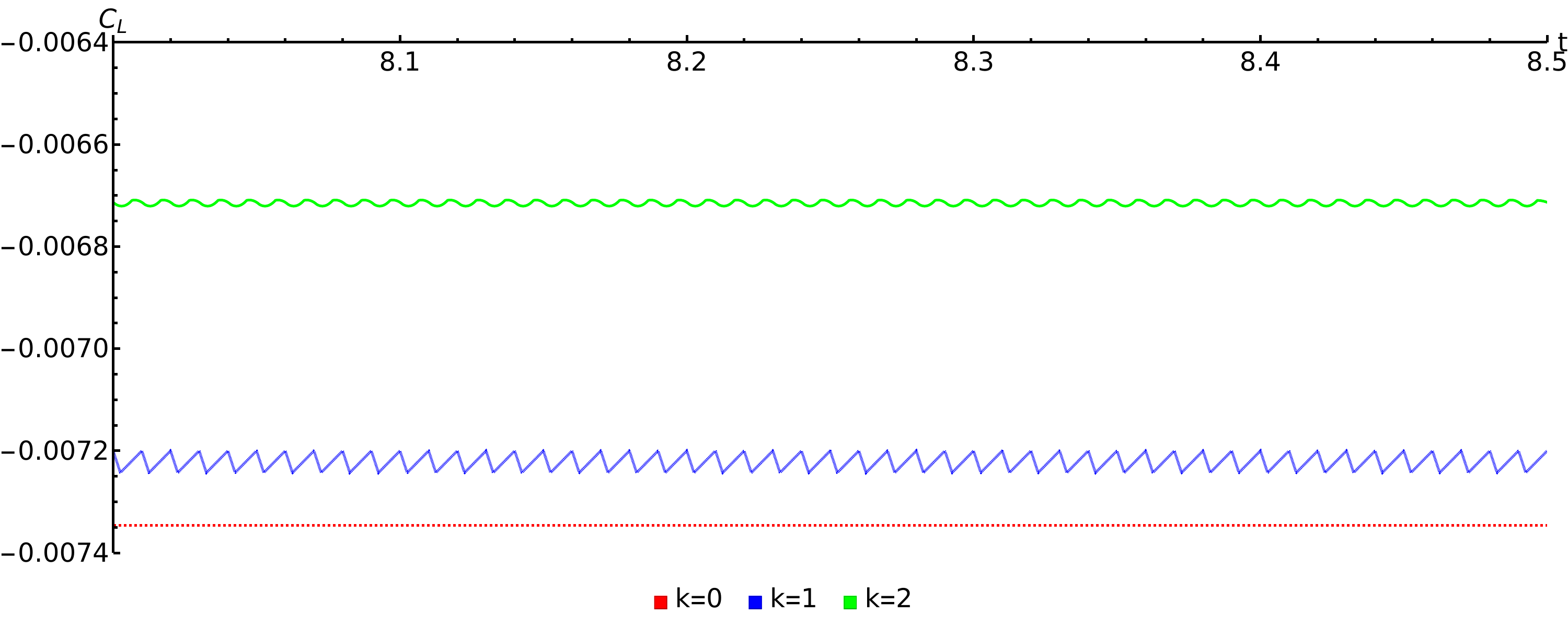}
	}
	\caption{Drag and lift coefficient for different time discretization schemes in the initial phase and in the fully developed flow regime.}%
	\label{fig:drag_lift_start_end}
\end{figure}

\begin{figure}[!htb]
	\centering
	\subcaptionbox{\centering Velocity profile at $t = T$ for $Re = 100$. \label{fig:plot_v_end_100}}
	[0.47\columnwidth]
	{\includegraphics[width=0.47\columnwidth,keepaspectratio]{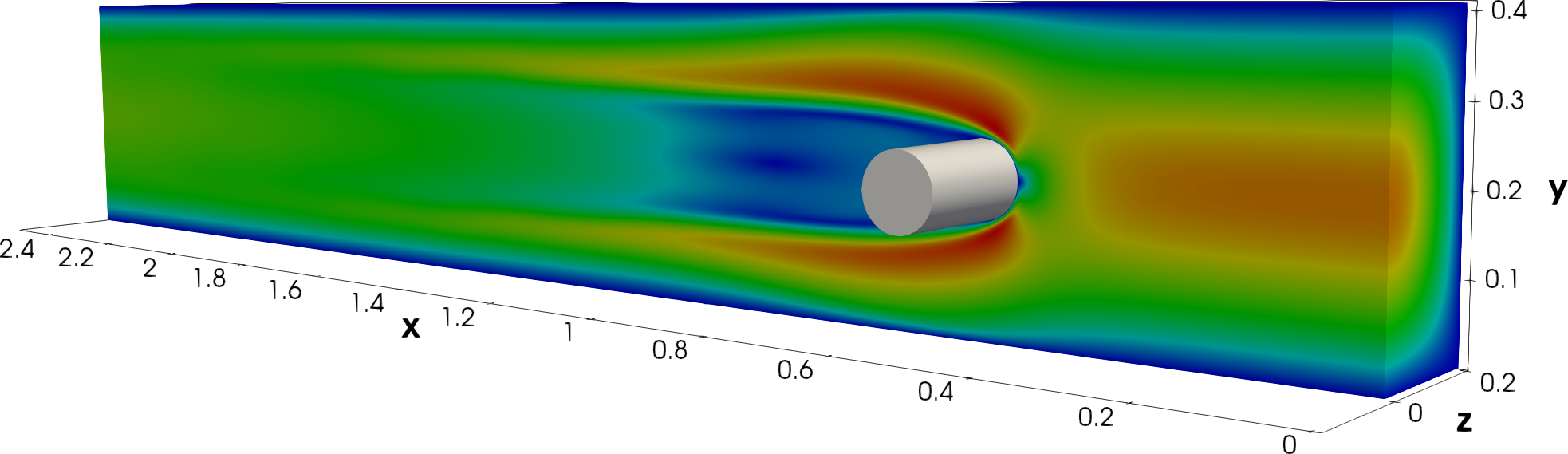}
	}
	\quad
	\subcaptionbox{\centering Velocity profile at $t = T$ for $Re = 133.\bar{3}$. \label{fig:plot_v_end_133}}
	[0.47\columnwidth]
	{\includegraphics[width=0.47\columnwidth,keepaspectratio]{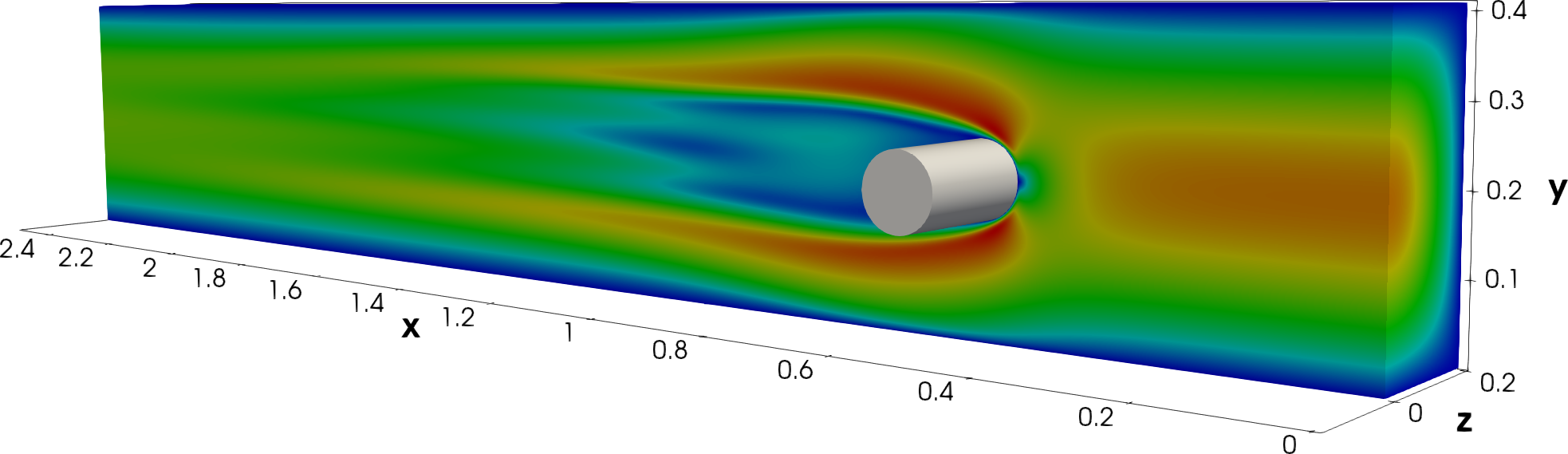}
	}
	\caption{Comparison of the velocity profile for flow at $Re = 100$ and $Re = 133.\bar{3}$ with $k=1$.}%
	\label{fig:velocity_100_133}
\end{figure}

Figure~\ref{fig:plot_v_end_100} illustrates the computed flow profile. For $Re = 100$, the flow profile approaches a stationary limit within the considered time interval. \Cref{fig:drag_lift_start_end} shows the computed drag and lift coefficient for the initial phase $t\in (0,0.5]$ and the final part $t\in (8.0,8.5]$ of the simulation for different values of the polynomial order $k$. In particular in the initial simulation interval of higher dynamics in the flow field a strong difference in the values for the drag and lift coefficient computed with different polynomial degrees for the time discretization can be observed; cf.~Fig.~\ref{fig:drag_start} and \ref{fig:lift_start}. Lower order temporal discretization schemes seem to damp the flow dynamics significantly, which in turn leads to different limits of the drag and lift coefficient at the end of the time interval of the simulation; cf.~Fig.~\ref{fig:drag_end} and \ref{fig:lift_end}.

For the three-dimensional case, the performance properties of the Newton-GMRES-GMG solver are summarized in Table~\ref{tab:3d_DFG_results}. Due to the transition of the flow profile into a steady state, the characteristic mean quantities of the Newton solver {$\bar{n}_{\text{Newton}}$} and of the GMRES method {$\bar{n}_{\text{GMRES}}$} are computed only until $t=t_s$ is reached, i.e.\ for $t\in (0,t_s]$. For $t>t_s$ the flow becomes almost stationary and nonlinear and linear iterations are no longer required. Again, the algebraic solver is highly efficient with performance properties being (almost) independent of the polynomial order $k$ of the time discretization. \Cref{tab:3d_DFG_results} also summarizes the computed maximum values of the drag and lift coefficient in the stationary regime.  

{
	\sisetup{scientific-notation = false,
		round-mode=places,
		round-precision=4,
		output-exponent-marker=\ensuremath{\mathrm{e}},
		table-figures-integer=1, 
		table-figures-decimal=3, 
		table-figures-exponent=1, 
		table-sign-mantissa = false, 
		table-sign-exponent = true, 
		table-number-alignment=center} 
	
	\begin{table}[h!t]
		\caption{Space-time degrees of freedom on each time interval $I_n$, timepoint $t_s$ when the flow becomes stationary, drag and lift coefficients as well as average number of Newton steps per time step and GMRES iterations per Newton step until $t=t_s$.}
		\centering
		\small
		\begin{tabular}{c@{\hskip 4ex} r@{\hskip 4ex} c@{\hskip 4ex} c@{\hskip 4ex}  c@{\hskip 4ex} c@{\hskip 4ex} c@{\hskip 4ex}}
			\toprule
			{$k$} & {DoFs$_{I_n}$} & {$t_{s}$} &{$c_{D_{max}}$} & {$c_{L_{max}}$} & {$\bar{n}_{\text{Newton}}$}   & {$\bar{n}_{\text{GMRES}}$} \\
			%
			\midrule
			{0} & 48438368 & {2.1} & \num{3.2877} & \num{-0.0072} & 1.43 & 2  \\
			{1} & 96876736 & {2.7} & \num{3.28588} & \num{-0.0070541} & 1.69 & 2  \\
			{2} & 145315104 & {3.6} &\num{3.28494} &  \num{-0.0065674} & 1.75 & 2 \\
			\bottomrule
		\end{tabular}
		\label{tab:3d_DFG_results}
	\end{table}
}

\begin{figure}[!htb]
	\centering
	\subcaptionbox{\centering Drag coefficient. \label{fig:drag_end_133}}
	[0.47\columnwidth]
	{\includegraphics[width=0.47\columnwidth,keepaspectratio]{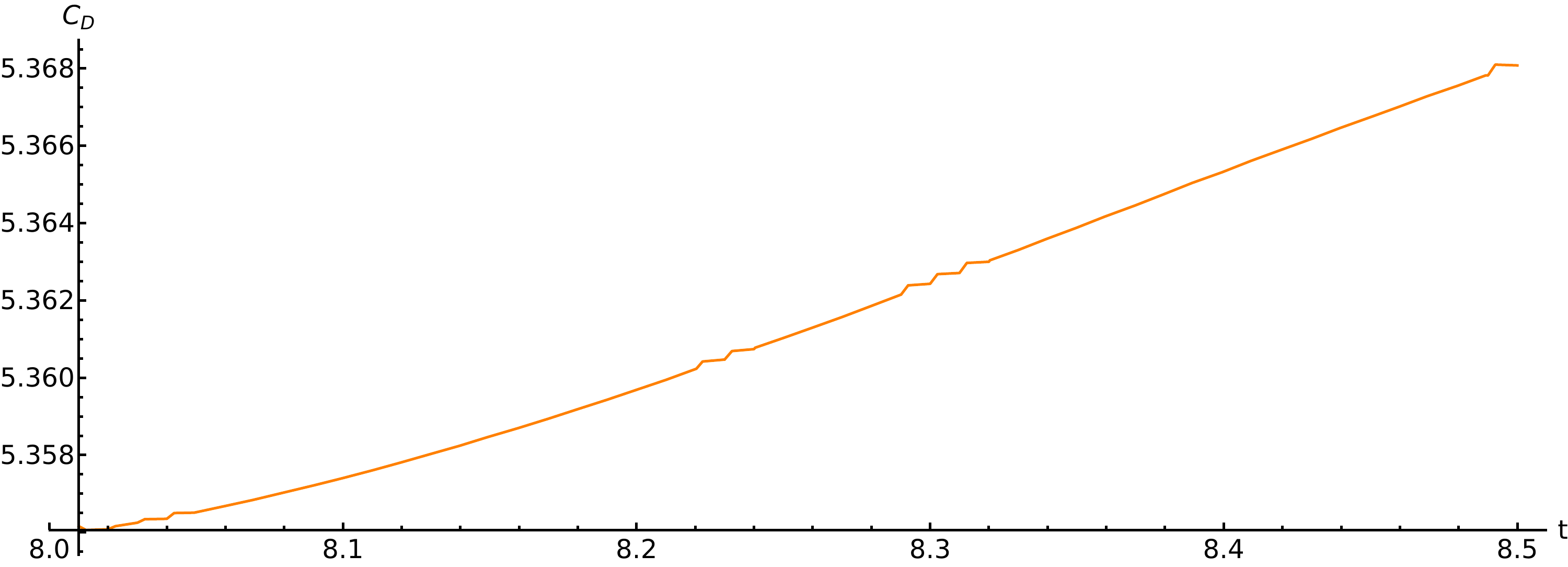}
	}
	\quad
	\subcaptionbox{\centering Lift coefficient. \label{fig:lift_end_133}}
	[0.47\columnwidth]
	{\includegraphics[width=0.47\columnwidth,keepaspectratio]{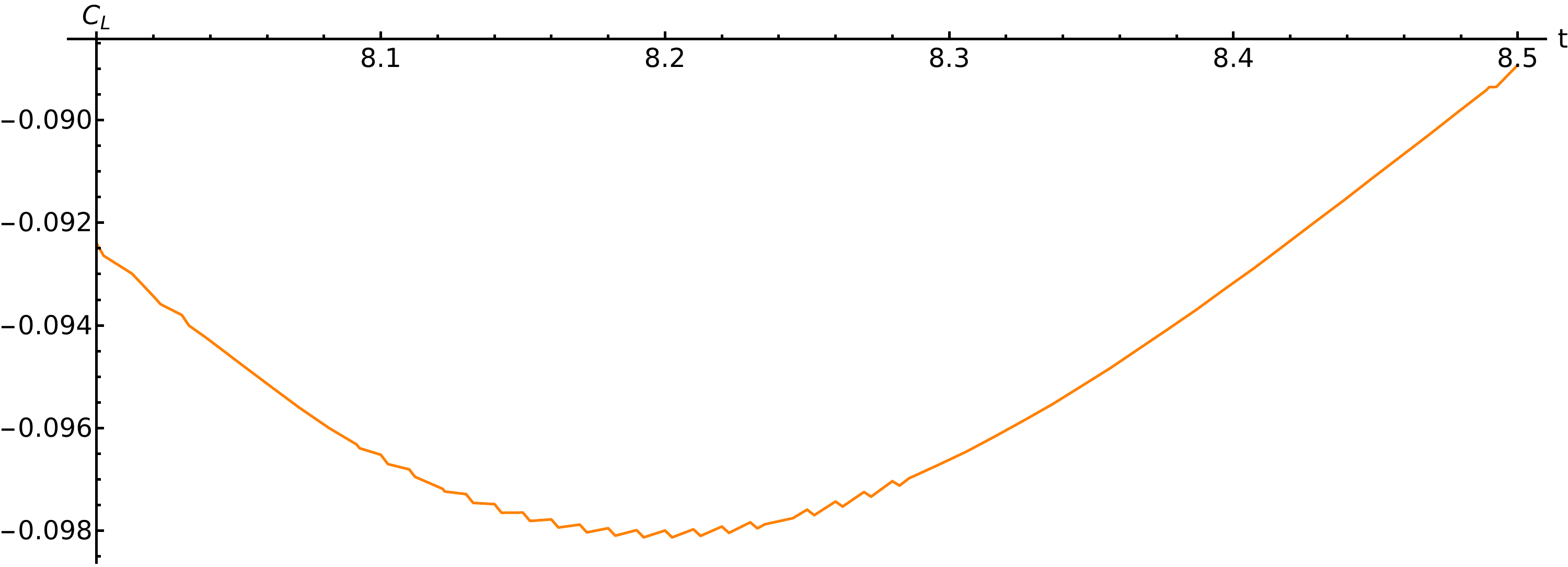}
	}
	\caption{Drag and lift coefficient for Reynolds number $Re =133.\bar{3}$, polynomial degree in time $k=1$ and time $t\in (8.0,8.5]$.}%
	\label{fig:drag_lift_end_133}
\end{figure}

Finally, in order to enhance the complexity of the flow problem and illustrate the transition to a nonstationary flow profile over the considered time interval, we replace the quantity $U_m = 2.25$ of the inflow profile of the benchmark setting (cf.~\cite{STDK96}) by $U_m$ = $3.0$, corresponding to Reynolds number $Re=133.\bar{3}$. Figure~\ref{fig:plot_v_end_133} shows the computed velocity profile and Figure~\ref{fig:drag_lift_end_133} the corresponding values of the drag and lift coefficient for $t\in (8.0,8.5]$ and polynomial degree $k=1$ of the time discretization. The plots show that the flow is nonstationary in the considered interval $(0,T]$, even though this is not obvious from Figure~\ref{fig:plot_v_end_133}. \Cref{tab:3d_DFG_results_comparison} compares the average number of Newton and GMRES iteration steps for the experiments with $Re=100$ and $Re=133.\bar{3}$. Since for $Re=100$ the flow gets stationary for simulation time of $t> 2.7$ (cf.\ Table~\ref{tab:3d_DFG_results}), the average number of Newton iterations throughout the entire time interval of the simulation is close to one. In contrast to this, for Reynolds number $Re = 133.\bar{3}$ the flow continues to remain nonstationary over the whole time interval and requires more Newton iterations in the time steps. Therefore, the average number ${\bar{n}_{\text{Newton}}}$ of Newton iterations is higher for $Re = 133.\bar{3}$ than for $Re = 100$, it is even higher than in the inital phase for $Re=100$; cf.\ Table~\ref{tab:3d_DFG_results}.

{
	\sisetup{scientific-notation = false,
		round-mode=places,
		round-precision=0,
		output-exponent-marker=\ensuremath{\mathrm{e}},
		table-figures-integer=1, 
		table-figures-decimal=3, 
		table-figures-exponent=1, 
		table-sign-mantissa = false, 
		table-sign-exponent = true, 
		table-number-alignment=center} 
	
	\begin{table}[h!t]
		\caption{Average number of Newton and GMRES iterations for $k=1$ and \num{96876736} DoFs on each subinterval $I_n$.}
		\centering
		\small
		\begin{tabular}{c@{\hskip 4ex}  c@{\hskip 4ex} c@{\hskip 4ex} c@{\hskip 4ex} c}
			\toprule
			{$Re$}  & {$\bar{n}_{\text{Newton}}$}  & Computed on  & {$\bar{n}_{\text{GMRES}}$} & Computed on\\
			%
			\midrule
			{100} & 1.00 & $(0,T]$ & 2  & $(0,t_s]$\\
			{133} & 1.79 & $(0,T]$ & 2  & $(0,T]$\\
			\bottomrule
		\end{tabular}
		\label{tab:3d_DFG_results_comparison}
	\end{table}
}

\section{Summary}

We presented an algebraic solver for space-time finite element approximations of the Navier--Stokes equations. For solving the block matrix systems of the Newton linearization, GMRES iterations with geometric multigrid preconditioning, using a local Vanka smoother, are applied. The efficiency of the solver and its robustness with respect to the (piecewise) polynomial degree of the time discretization was studied and demonstrated for the 2d and 3d benchmark problem of flow around a cylinder.

\end{document}